\documentclass[11pt]{article}
\usepackage{euler}
\usepackage{ams}
\usepackage{xypic}
\newtheorem{la}{Lemma}
\newtheorem{sz}{Proposition}
\newtheorem{tm}{Theorem}

\input xy
\xyoption{all}

\newcommand{\Sp} [0] {{\rm Tr}}

\begin{document}
\begin{center}{\Large{\bf On the Zeros of Fermat Quotients and Mirimanoff Polynomials}}\\
\vspace{.5cm}
{Bj\"orn Grohmann\footnote{email: nn@mhorg.de}}\\
{April 2006}\\
\end{center}
\vspace{.5cm}
\begin{abstract}
In this article an upper bound for the first consecutive zeros of the Fermat quotient is given in terms of the zeros of a Mirimanoff
polynomial. This bound is obtained by investigating a relation between these polynomials and the factor system of Gauss sums of a 
certain number field.\\\\
{\bf Keywords:} Mirimanoff polynomials, Fermat quotients, Gauss sums.\\ 
\end{abstract}
\section{Introduction}
Let $p$ be an odd prime and $k$ an integer, with $k\not\equiv 0\mbox{ mod }p$. The {\em Fermat quotient} $q_p(k)$ is defined to be
the smallest positive integer, which satisfies the equation
\begin{equation}
k^{p-1}\equiv 1+q_p(k)p\mbox{ mod }p^2.
\end{equation}
During the last century, many properties of this quotient have been investigated, mostly because of its connection to the
first case of Fermat's Last Theorem: In 1909 Wieferich was able to show that if there exists an odd prime $p$ and nonzero integers $x,y,z$, with
$p\!\!\not|\,\, xyz$ and $x^p+y^p+z^p=0$, then 
\begin{equation}\label{qp2}
q_p(2)=0.
\end{equation}
Until today, only two primes are known to satisfy equation (\ref{qp2}): the prime $1093$ (found by W. Meissner in 1913) 
and $3511$ (N. Beeger, 1921/22).\\\\
One year after Wieferich, Mirimanoff stated that in the above case, it also holds that $q_p(3)=0$, and it is still an open question 
(independent of {\rm FLT}), whether there exists a prime $p$, such that
\begin{equation}
q_p(2)=q_p(3)=0.
\end{equation}
The interested reader should consult \cite{ribenboim79} for more on this fascinating story.\\\\
In this article we will derive a new upper bound for the integer
\begin{equation}
\kappa_p:=\min\{n>0\,|\,q_p(n)\not=0\},
\end{equation}
where $p$ denotes an odd prime. As can be easily seen from the definiton of $q_p(\cdot)$, it holds that
\begin{equation}\label{qplog}
q_p(ab)\equiv q_p(a)+q_p(b)\mbox{ mod }p,
\end{equation}
so the integer $\kappa_p$ is prime too.\\\\
The bound for $\kappa_p$ will be in terms of the number of zeros of the {\em Mirimanoff polynomial}
\begin{equation}
\gamma_p(t):=\sum_{j=1}^{p-1}\frac{t^j}{j}
\end{equation}
modulo $p$. This polynomial is closely related to the Fermat quotient, since
\begin{equation}\label{mpolyr1}
\gamma_p(t)
\equiv \frac{1-t^p-(1-t)^p}{p}
\equiv (t-1)q_p(t-1)-tq_p(t)\mbox{ mod }p,
\end{equation}
and therefore 
\begin{equation}\label{kpgp}
\kappa_p=\min\{n>0\,|\,\gamma_p(n)\not\equiv 0\mbox{ mod }p\}.
\end{equation}
In the sequel we will denote the number of zeros of $\gamma_p$ modulo $p$ by $\eta_{0,p}$ and in general, for $0\leq a<p$:
\begin{equation}
\eta_{a,p}:=|\{c\mid 0\leq c <p, \gamma_p(c)\equiv a\mbox{ mod } p\}|.
\end{equation}
Here, we will prove the following bound:
\vspace{0.5cm}
\begin{tm}\label{kpbound}
\begin{equation}
\kappa_p\in O(\sqrt{\eta_{0,p}}).
\end{equation}
\end{tm}
\vspace{0.5cm}
We will start by recalling some basic facts about Mirimanoff polynomials and deriving some trivial bounds for $\kappa_p$.
After that, we will state a connection between the factor systems of Gauss sums of a certain number field with these polynomials.
This will lead to some useful relations that will enable us to prove the stated result.
\section{Basic Properties}
In the following, when we talk about the {\em zeros} of polynomial $\gamma_p$, we will always mean the zeros modulo $p$.\\\\
The first thing to note is that $0$ and $1$ are zeros of $\gamma_p$ and we will call them the trivial zeros. All the other zeros,
if they exist, have multiplicity two, since $\gamma_p^{\prime}(t)=1+t+\dots+t^{p-2}$. From this and (\ref{kpgp}) it is immediate that
\begin{equation}\label{kapab2}
\kappa_p\leq \frac{p+1}{2}.
\end{equation}
Looking at equation (\ref{mpolyr1}), we see that the Mirimanoff polynomial satisfies
\begin{equation}\label{gpema}
\gamma_p(a)\equiv \gamma_p(1-a)\mbox{ mod }p,
\end{equation}
and, for $a\not\equiv 0\mbox{ mod }p$, 
\begin{equation}\label{gpeda}
\gamma_p(a)\equiv -a\gamma_p(\frac{1}{a})\mbox{ mod }p.
\end{equation}
From this it follows that a nontrivial zero $z$ of $\gamma_p$ always leads to the zeros:
\begin{equation}
\xymatrix@R-13pt@C-5pt{&\frac{z}{1}\ar @{.} [dl]\ar @{-}[dr]&\\
\frac{1}{z}\ar @{-} [d]&&\frac{1-z}{1} \ar @{.} [d]\\
\frac{z-1}{z}\ar @{.}[dr]&&\frac{1}{1-z}\ar @{-}[dl]\\
&\frac{z}{z-1}&
}
\end{equation}
There are two exceptional cases. In case of $\gamma_p(2)\equiv 0\mbox{ mod }p$ we have
\begin{equation}
\xymatrix@R-13pt@C-5pt{&2\ar @{.} [dl]\ar @{-}[dr]&\\
\frac{1}{2}\ar@(ul,dl) @{-}&&-1 \ar@(ur,dr) @{.}
}
\end{equation}
and in case of $p\equiv 1 \mbox{ mod }3$ it is easy to see that the zeros $\alpha_6,\alpha_6^5$ of the polynomial $t^2-t+1\in \F_p[t]$
are aswell zeros of $\gamma_p$, since
\begin{equation}
\xymatrix@R-13pt@C-5pt{
\alpha_6\ar@/_/ @{.} [rr]\ar@/^/ @{-}[rr]&&\alpha_6^5
}.
\end{equation}
By what has been said so far, it is clear that only half of the zeros of $\gamma_p$ contribute to an upper bound of $\kappa_p$, which means
that
\begin{equation}\label{kaab3}
\kappa_p\leq \left\lfloor \frac{p+5}{4} \right\rfloor.
\end{equation}
\section{Relations}
The Mirimanoff polynomials $\gamma_p$ satisfy a number of remarkable relations. To give only one example: If $z$ is a nontrivial zero and
if in addition
\begin{equation}
\gamma_p(z)\equiv \gamma_p(z+1)\equiv 0 \mbox{ mod }p,
\end{equation}
then it is also true that
\begin{equation}
\gamma_p(z^2)\equiv 0 \mbox{ mod }p.
\end{equation}
The aim of this section is to prove the following relation.
\vspace{0.5cm}
\begin{sz}\label{mpmr}
For integers $a,b$, with $a,b,ab\not\equiv 0,1\mbox{ {\rm mod} }p$ we have
\begin{eqnarray}\label{fundrelg}
(1-ab)\gamma_p(\frac{1-b}{1-ab})+(1-b)\gamma_p(a)\equiv 
(1-ab)\gamma_p(\frac{1-a}{1-ab})+(1-a)\gamma_p(b)
\mbox{ {\rm mod} }p.
\end{eqnarray} 
\end{sz}
\vspace{0.5cm}
Here, we obtain the above example by setting $a:=z$ and $b:=1/(z+1)$.\\\\
To prove the proposition, we introduce the abelian number field $K$, which is the real extension of the field of rationals $\Q$ of degree
$(K:\Q)=p$ whose discriminant is a power of $p$. It may be described as the subfield $\Q(\Gamma)$ of the cyclotomic 
extension $\Q(\zeta_{p^2})/\Q$, where
$\zeta_{p^2}:=e^{2\pi i/p^2}$ is a primitive $p^2$-th root of unity, and 
\begin{equation}\label{defggamma}
\Gamma:=\Sp_{\Q(\zeta_{p^2})/K}(\zeta_{p^2}).
\end{equation} 
Hence, all nontrivial characters $\chi$ of $K$ have conductor $f_\chi=p^2$ and are even ($\chi(-1)=1$). As usual, we will denote
the Gauss sum of a character $\chi$ of $K$ by $\tau(\chi)$ and the elements of the factor system of these sums by
\begin{equation}
\omega(\chi,\chi^{\prime}):=\frac{\tau(\chi)\tau(\chi^{\prime})}
{\tau(\chi\chi^{\prime})},
\end{equation}
for characters $\chi, \chi^\prime$ of $K$.
For an introduction to characters of abelian number fields and Gauss sums we refer the reader to \cite{washington97},\cite{berndt98}.\\\\
We now fix a generator $\chi$ of the (cyclic) character group of $K$ and assume that $\chi$ is normalized
\begin{equation}
\chi(1-p)=\zeta_p.
\end{equation}
There is the following connection between the elements of the factor system and the Mirimonoff polynomial $\gamma_p$.
\vspace{0.5cm}
\begin{sz}\label{wgsz} 
For $k,l,k+l \not\equiv 0\mbox{ {\rm mod} }p$ it holds that
\begin{equation}
\omega(\chi^k,\chi^l)=p\zeta_p^{-(k+l)\gamma_p(\frac{k}{k+l})}.
\end{equation}
\end{sz}
\vspace{0.5cm}
{\bf Proof.} We start by evaluating the Gauss sum $\tau(\chi)$. Since $f_\chi=p^2$ and $\chi(1-p)=\zeta_p$, Theorem 1.6.2 from \cite{berndt98} gives
\begin{equation}
\tau(\chi)=p\zeta_{p^2}.
\end{equation}
Now, for an integer $k$, with $(k,p)=1$, denote by $\sigma_k$ the automorphism of the Galois group $G(\Q(\zeta_{p^2})/\Q)$ definied
by $\zeta_{p^2}^{\sigma_k}:=\zeta_{p^2}^k$. It is then well known that (cf. \cite{mydiss})
$\tau(\chi)^{\sigma_k}=\bar\chi^k(k)\tau(\chi^k)$,
which leads to
\begin{equation}
\tau(\chi^k)=\chi^k(k)p\zeta_{p^2}^k.
\end{equation}
The elements of the factor system may therefore be written as
\begin{equation}
\omega(\chi^k,\chi^l)=p\frac{\chi^k(k)\chi^l(l)}{\chi^{k+l}(k+l)}.
\end{equation} 
Since the character is normalized it follows that
\begin{equation}
\chi(k)=\zeta_p^{q_p(k)},
\end{equation}
so by (\ref{qplog}) and (\ref{mpolyr1}) we obtain the stated result.\hfill$\Box$\\\\
{\bf Proof of Proposition \ref{mpmr}.} As is easily seen, the elements of the factor system satisfy the 2-cocycle condition
\begin{equation}
\omega(\chi,\chi^{\prime}\chi^{\prime\prime})\omega(\chi^{\prime},\chi^{\prime\prime})=
\omega(\chi\chi^{\prime},\chi^{\prime\prime})\omega(\chi,\chi^{\prime}),
\end{equation}
for all characters $\chi, \chi^\prime, \chi^{\prime\prime}$ of $K$. Using Proposition \ref{wgsz}, this leads for
integers $k,l,j$, with $k,l,j,k+l,l+j,k+l+j\not\equiv 0\mbox{ mod }p$, to
\begin{eqnarray*}
(k+l+j)\gamma_p(\frac{k}{k+l+j})+(l+j)\gamma_p(\frac{l}{l+j})\equiv
(k+l+j)\gamma_p(\frac{k+l}{k+l+j})+(k+l)\gamma_p(\frac{k}{k+l}).
\end{eqnarray*}
We now obtain the relation of Proposition \ref{mpmr} by substituting $1-a:=l/(l+j)$ and $b:=k/(k+l)$.\hfill$\Box$\\\\
\section{Proof of the Main Theorem}
As a consequence of Proposition \ref{mpmr} we get the following statement.
\vspace{0.5cm}
\begin{sz}\label{gpuv}
For $1\leq u,v <\kappa_p$ it holds that
\begin{equation}
\gamma_p(\frac{u}{v})\equiv 0 \mbox{ {\rm mod} }p.
\end{equation}
\end{sz}
\vspace{0.5cm}
{\bf Proof.}  Deviding equation (\ref{fundrelg}) by $(1-a)(1-b)$ and substituting $e:=1/(1-a)$ and $d:=1/(1-b)$ leads,
for $e,d\not\equiv 0,1\mbox{ mod }p$ and $e+d\not\equiv 1\mbox{ mod }p$, to
\begin{equation}
(e+d-1)\left(\gamma_p(\frac{d-1}{e+d-1})-\gamma_p(\frac{d}{e+d-1})\right)\equiv \gamma_p(e)-\gamma_p(d) \mbox{ mod }p.
\end{equation}
The statement now easily follows by induction, taking into account the relations (\ref{gpema}), (\ref{gpeda}) and 
(\ref{kaab3}).\hfill$\Box$\\\\
{\bf Proof of Theorem \ref{kpbound}.} To prove the bound of  Theorem \ref{kpbound} we first note that, for any positive integer $q$,
the integer
\begin{equation}
s_q:=|\{(u,v)\mid 1\leq u,v \leq q,\, {\rm gcd}(u,v)=1\}|
\end{equation}
satisfies
\begin{equation}\label{sqbound}
s_q\geq\sum_{k=1}^q q-q\left(\sum_{l|k}\frac{1}{l}\right)\geq q^2\left(1-\sum_{l\leq q}\frac{1}{l^2}\right)\geq 
q^2\left(2-\frac{\pi^2}{6}\right).
\end{equation}
We now claim that for large $p$
\begin{equation}
\kappa_p\leq \lfloor\sqrt p \rfloor.
\end{equation}
It then follows from Propostition \ref{gpuv} and (\ref{sqbound}) that there exist $\Omega(\kappa_p^2)$ zeros of $\gamma_p$ and therefore
the bound of the Theorem holds.\\\\
To prove the claim, we assume that $\kappa_p>\lfloor\sqrt p \rfloor$. This leads to more than
\begin{equation}
\left(\lfloor\sqrt p\rfloor\right)^2\left(2-\frac{\pi^2}{6}\right)>\frac{p}{3}
\end{equation}
zeros of the form $a/b$, with $(a,b)=1$ and $1\leq a,b\leq \lfloor\sqrt p \rfloor$. Now, if $a/b$ is a zero of $\gamma_p$,
so is $1-a/b$. We will need the following Lemma.
\begin{la}
Let $p$ be an odd prime and $a,b$ integers, with $1\leq a,b \leq \lfloor\sqrt p \rfloor$ and $a>b$. Then there do {\underline{not}} exist
integers
$c,d$, with $1\leq c,d \leq \lfloor\sqrt p \rfloor$, such that:
\begin{equation}
1-\frac{a}{b}\equiv \frac{c}{d}\,\,{\rm mod}\,p.
\end{equation}
\end{la}
It follows that for large $p$ and under the assumption $\kappa_p>\lfloor\sqrt p \rfloor$, we have produced
\begin{equation}
\bigg\lceil\frac{p}{3}\bigg\rceil+\bigg\lceil\frac{1}{2}\left(\bigg\lceil\frac{p}{3}\bigg\rceil -1\right)\bigg\rceil=\frac{p+1}{2}
\end{equation}
zeros of $\gamma_p$, which is impossible, since we did not count $0$ and all nontrival zeros have multiplicity two. This means that
$\kappa_p\leq \lfloor\sqrt p \rfloor$ and therefore the Theorem follows.\\\\
{\bf Proof of the Lemma.} From the equation
\begin{equation}
bd-ad\equiv bc\,\,{\rm mod}\,p
\end{equation}
and $a>b$ it follows that there exist an integer $k>0$, such that
\begin{equation}
bd-ad+kp=bc.
\end{equation}
Now, $b$ devides $kp-ad$, say $kp-ad=sb$, with $s>0$ and in particular $d+s=c$, so $c>d$ and
\begin{equation}
1\leq s\leq\lfloor\sqrt p\rfloor.
\end{equation}
This leads to
\begin{equation}
kp=sb+ad<2p,
\end{equation}
so $k=1$ and therefore $bd-ad+p=bc$, resp. $d=(p-bc)/(a-b)$. Since $c>d$ it follows that $ac>p$, which contradicts the assumptions of
the Lemma.\hfill$\Box$\\\\ 

\end{document}